\documentclass[11pt]{amsart}
\usepackage{amstext}
\usepackage{amsfonts}
\usepackage{amssymb}
\usepackage{amsbsy}
\usepackage{latexsym}
\usepackage[T1]{fontenc}
\usepackage{xy}
\usepackage{hhline}
\xyoption{all}
\newcommand{\vtx}[1]{*+[o][F-]{\scriptscriptstyle #1}} 
\newcounter{num}[section] %

\newenvironment{theo}
{\refstepcounter{num}%
\bigskip\noindent{\bf Theorem~\arabic{section}.\arabic{num}. }\it}
{\smallskip}

\newcommand{\mylabel}[1]{}

\newenvironment{eq}{\begin{equation}}{\end{equation}}

\newcommand{\Ref}[1]{(\ref{#1})}

\newcommand{\si}{\sigma}
\newcommand{\al}{\alpha}

\newcommand{\la}{\lambda}

\newcommand{\LA}{\langle}
\newcommand{\RA}{\rangle}
\newcommand{\ov}[1]{\overline{#1}}
\newcommand{\un}[1]{{\underline{#1}} }

\newcommand{\tr}{\mathop{\rm tr}}

\newcommand{\Char}{\mathop{\rm char}}

\newcommand{\Hom}{{\mathop{\rm{Hom }}}}

\newcommand{\alg}{{\rm alg}}
\newcommand{\algA}{\mathcal{A}}    

\newcommand{\FF}{{\mathbb{F}}}   
\newcommand{\NN}{{\mathbb{N}}}

\newcommand{\Q}{\mathcal{Q}}    
\newcommand{\n}{\boldsymbol{n}} 
 
\renewcommand{\i}{\boldsymbol{i}} 
\newcommand{\Sp}{S\!p}
\newcommand{\D}{{D}}

\newcommand{\PhiLarge}[1]{\widetilde{\Phi}_{#1}}  

\newcommand{\Y}{\LA Y\RA}
\newcommand{\YQ}{\LA Y_{\Q}\RA}
\newcommand{\EY}{\LA \widetilde{Y}\RA}

\newcommand{\AlgLargeY}{\si\EY}  
\newcommand{\KLargeY}[1]{\widetilde{K}'_{#1}}  
\newcommand{\PhiLargeY}[1]{\widetilde{\Phi}'_{#1}}  
\newcommand{\piLargeY}{\widetilde{\pi}'}  

\newcommand{\AlgSmallY}[1]{\si\Y_{#1}}  
\newcommand{\KSmallY}[1]{K'_{#1}}  
\newcommand{\PhiSmallY}[1]{\Phi'_{#1}}  
\newcommand{\piSmallY}[1]{\pi'_{#1}}  

\newcommand{\PhiAbsY}[1]{\widehat{\Phi}'_{#1}}   

\newcommand{\AlgLargeQ}{\si\LA\widetilde{Q}\RA}  
\newcommand{\KLargeQ}[1]{\widetilde{K}_{#1}^{\Q}}  
\newcommand{\PhiLargeQ}[1]{\widetilde{\Phi}_{#1}}  
\newcommand{\piLargeQ}{\widetilde{\pi}}  
\newcommand{\EQ}{\LA \widetilde{Q}\RA}

\newcommand{\AlgSmallQ}[1]{\si\LA \Q\RA_{#1}} 
\newcommand{\KSmallQ}[1]{K^{\Q}_{#1}}  
\newcommand{\PhiSmallQ}[1]{\Phi_{#1}}  
\newcommand{\piSmallQ}[1]{\pi_{#1}}  

\newcommand{\PhiAbsQ}[1]{\widehat{\Phi}_{#1}}   

\begin{document}
\title[Identities for mixed representations of quivers are finitely based]{Identities for mixed representations of quiver are finitely based}
 \author{Artem Lopatin}
\thanks{This research was supported by RFBR 16-31-60111 (mol\_a\_dk).}
\address{Artem Lopatin\\
Sobolev Institute of Mathematics, Omsk Branch, SB RAS, Omsk, Russia} 
\email{artem\_lopatin@yahoo.com}

\begin{abstract} 
The classical construction of representations of quivers enables us to consider linear maps between several vector spaces. The mixed representations of quivers helps us to work with linear maps as well as bilinear forms on several vector spaces. We construct a finite generating set for the T-ideal of identities of the polynomial invariants of mixed representations of a quiver. 



\end{abstract}

\maketitle

\section{Introduction}\label{section_intro}

All vector spaces, algebras, and modules are over an infinite field $\FF$ of characteristic $\Char{\FF}=p\geq0$.  By an algebra we always mean an associative algebra. We write $\FF^{n\times m}$ for the set of all $n\times m$ matrices over $\FF$. Given an $n\times n$ matrix $A$ we write $A_{ij}$ for the $(i,j)^{\rm th}$ entry of $A$. We set $\NN_0=\{0,1,2,\ldots\}$ and $|\un{t}|=t_1+\cdots+t_u$ for $\un{t}\in\NN_0^u$.

\subsection{Invariants of bilinear forms} 

Consider bilinear forms $(\cdot,\cdot)_1,\ldots,(\cdot,\cdot)_r$ on the vector space  $V=\FF^n$ , which are defined by $n\times n$ matrices $B_1,\ldots,B_r$, respectively, with the respect of some fixed base of $V$. Similarly, we consider bilinear forms  $\LA\cdot,\cdot\RA_1,\ldots,\LA\cdot,\cdot\RA_s$ on the dual space $V^{\ast}$, which are defined by $n\times n$ matrices $C_1,\ldots,C_s$, respectively, with the respect of the dual bases of $V^{\ast}$.  Then $G=GL(n)$ acts on the space 
$$H=H_{r,s}=\bigoplus_{k=1}^r \FF^{n\times n}\oplus \bigoplus_{l=1}^s \FF^{n\times n}$$ 
of the above mentioned bilinear forms as base change: 
$$g\cdot (B_1,\ldots,B_r,C_1,\ldots,C_s)=(g B_1 g^{T},\ldots,g B_r g^{T},g^{-T}C_1 g^{-1},\ldots,g^{-T}C_s g^{-1}),$$
where $g^{-T}$ stands for $(g^T)^{-1}$. This action induces the natural action of $GL(n)$ on the coordinate ring of the affine variety $H$
$$\FF[H]=\FF[x_{ij}(k),\,y_{ij}(l) \,|\,1\leq i,j\leq n,\,1\leq k\leq r, \,1\leq l\leq s]$$
by the formula: $(g\cdot f)(h)=f(g\cdot h)$ for $g\in G$, $f\in \FF[H]$, and $h\in H$, where we consider $f$ as a polynomial function from $H$ to $\FF$. To write down the explicit formula for this action, denote the {\it generic $n\times n$ matrices} by $X_k=(x_{ij}(k))$ and $Y_l=(y_{ij}(l))$. Then $g\cdot x_{ij}(k) = (g^{-1} X_k g^{-T})_{ij}$ and $g\cdot y_{ij}(l) = (g^{T} Y_l g)_{ij}$.
The set of elements of $\FF[H]$ that are stable with the respect to the action of $GL(n)$ is called the {\it algebra of invariants} of bilinear forms on $V$ and $V^{\ast}$ 
$$\FF[H]^{GL(n)}=\{f\in \FF[H]\,|\, g\cdot f=f\text{ for every }g\in GL(n)\}.$$
In other words, $\FF[H]^{GL(n)}=\{f\in \FF[H]\,|\, f(g\cdot h)=f(h)\text{ for every }g\in GL(n),\; h\in H\}$. It was shown in~\cite{ZubkovI} that $\FF[H]^{GL(n)}$ is generated by $\si_t(Z_1\cdots Z_m)$, where $1\leq t\leq n$ and $Z_i$ is one of the following products:
\begin{eq}\label{eq1} \mylabel{eq1}
X_k Y_l,\,X_k^T Y_l,\,X_k Y_l^T,\,X_k^T Y_l^T\; (1\leq k\leq r,\,1\leq l\leq s)
\end{eq}
and $\si_t(A)$ is the $t^{\rm th}$ coefficient of the characteristic polynomial of an $n\times n$ matrix $A$:
$$\det(\la E - A) = \la^n - \si_1(A) \la^{n-1} + \cdots + (-1)^n \si_n(A),$$
i.e. $\si_1(A)=\tr(A)$ and $\si_n(A)=\det(A)$. Moreover, we can assume that the product $Z_1\cdots Z_m$ is {\it primitive}, i.e., is not equal to the power of a shorter product of elements~\Ref{eq1}. Relations between these generators modulo free relations were described in~\cite{ZubkovII} and complete description of relations was obtained in Theorem~6.2 of~\cite{Lopatin_Ofree}. The  algebra of invariants of bilinear forms on $V$ and $V^{\ast}$ is a partial case of more general construction of invariants of mixed representations of quivers (see Example~\ref{ex1} below). 


\subsection{Invariants of mixed representations of quivers}

A {\it quiver} $\Q=(\Q_0,\Q_1)$ is a finite oriented graph, where $\Q_0$ ($\Q_1$, respectively) stands for the set of
vertices (the set of arrows, respectively). For an arrow $a$, denote by $a'$ its head
and by $a''$ its tail. 
We say that $a=a_1\cdots a_r$ is a {\it path} in $\Q$ (where $a_1,\ldots,a_r\in
\Q_1$), if $a_1''=a_2',\ldots,a_{r-1}''=a_r'$. 
The head of the path $a$ is $a'=a_1'$ and the tail is $a''=a_r''$.  A path
$a$ is called {\it closed} if $a'=a''$.

Given a {\it dimension vector} $\n=(\n_v\,|\,v\in\Q_0)$, we consider 
\begin{enumerate}
\item[$\bullet$] the space $H=\sum_{a\in\Q_1} \FF^{\n_{a'}\times \n_{a''}}\simeq \sum_{a\in\Q_1} \Hom(\FF^{\n_{a''}},\FF^{\n_{a'}})$;

\item[$\bullet$] the coordinate ring $R=\FF[x_{ij}^a\,|\,a\in\Q_1,\,1\leq i\leq \n_{a'},\,1\leq j\leq \n_{a''}]$ of $H$;

\item[$\bullet$] the $\n_{a'}\times \n_{a''}$ {\it generic} matrix $X_a=(x_{ij}^a)$ for every $a\in\Q_1$; 

\item[$\bullet$] the group $GL(\n)=\sum_{v\in\Q_0} GL(\n_v)$, acting on $H$ as the base change, i.e., 
$$g\cdot (h_a)=(g_{a'} h_a g^{-1}_{a''})$$
for $g=(g_v)\in GL(\n)$ and $(h_a)\in H$; this action induces the action of $GL(\n)$ on $R$. 
\end{enumerate}
The action of $GL(\n)$ on $H$ naturally induces the action on $R$. The set of stable elements with the respect of this action is called the algebra of {\it invariants} of representations of $\Q$
$$I(\Q,\n)=R^{GL(\n)}.$$
Given a path $a=a_1\cdots a_r$ with $a_i\in\Q_1$, we write $X_a$ for $X_{a_1}\cdots X_{a_r}$. Donkin~\cite{Donkin94} proved that $I(\Q,\n)$ is the subalgebra of $R$ generated by $\si_t(X_a)$, where $a$ is a closed path in $\Q$ and $1\leq t\leq \n_{a'}$. Moreover, we can assume that $a$ is {\it primitive}, i.e., is not equal to the power of a shorter closed path in $\Q$.

Let $\i:\Q_0\to\Q_0$ be an involution, i.e., $\i^2$ is the identical map, satisfying $\i(v)\neq v$ and $\n_{\i(v)}=\n_v$ for every vertex $v\in\Q_0$. Define 
\begin{enumerate}
\item[$\bullet$] the group $GL(\n,\i)\subset GL(\n)$ by 
$(g_v)\in GL(\n,\i)$ if and only if $g_v g_{\i(v)}^T=E$ for all $v$;

\item[$\bullet$] the {\it double} quiver $\Q^{\D}$ by $\Q_0^{\D}=\Q_0$ and  $\Q_1^{\D}=\Q_1\coprod \{a^T\,|\,a\in \Q_1\}$, 
where $(a^T)'=\i({a''})$, $(a^T)''=\i({a'})$ for all $a\in\Q_1$.  
\end{enumerate}
We set $X_{a^T}=X_a^T$ for all $a\in\Q_1$. Zubkov~\cite{ZubkovI} showed that the algebra of {\it invariants} of {\it mixed} representations of $\Q$
$$I(\Q,\n,\i)=R^{GL(\n,\i)}$$
is the subalgebra of $R$ generated by $\si_t(X_a)$, where $a$ is a closed path in $\Q^D$ and $1\leq t\leq \n_{a'}$. As above, we can assume that $a$ is primitive.  

\begin{example}\label{ex1}
Let $\Q$ be the following quiver
$$\vcenter
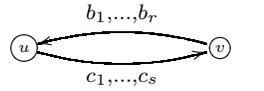
{
\xymatrix@C=1cm@R=1cm{ %
\vtx{u}\ar@/^/@{<-}[rr]^{b_1,\ldots,b_r} \ar@/_/@{->}[rr]_{c_1,\ldots,c_s}&&\vtx{v}\\
}} \quad,
$$
where there are $r$ arrows from $v$ to $u$ and $s$ arrows in the opposite direction,  
$\i(u)=v$, and $\n=(n,n)$. Then the algebra of invariants $\FF[H_{r,s}]^{GL(n)}$ of bilinear forms is isomorphic to $I(\Q,\n,\i)$. 
\end{example}

\subsection{The known results}

Assume that $G$ is one of the following groups: $GL(n)$, $O(n)$, $\Sp(n)$, where in case of $O(n)$ we assume $p\neq2$ and in case $\Sp(n)$ we have that $n$ is even. The algebra of matrix $G$-invariants $R^G$ in known to be generated by 
\begin{enumerate}
\item[$\bullet$] $\si_t(X_{i_1}\cdots X_{i_r})$ for $1\leq t\leq n$ in case $G=GL(n)$; 

\item[$\bullet$] $\si_t(A_{i_1}\cdots A_{i_r})$ for $1\leq t\leq n$, where $A_i$ is $X_i$ or $X_i^T$, in case $G=O(n)$ and $p\neq 2$; 

\item[$\bullet$] $\si_t(B_{i_1}\cdots B_{i_r})$ for $1\leq t\leq n$, where $B_i$ is $X_i$ or $X_i^T$, in case $G=\Sp(n)$ and $n$ is even. 
\end{enumerate}
Here $X_i$ stands for $n\times n$ generic matrix defined in Section~1.1. Note thatiIn case $p=0$ it is enough to consider traces instead of $\si_t$. In case of characteristic zero these results were obtained in~\cite{Sibirskii_1968},~\cite{Procesi76} and relations between these generators were independently computed in~\cite{Razmyslov74},~\cite{Procesi76}. Generators in case of positive characteristic were obtained in~\cite{Donkin92a},~\cite{Zubkov99} and relations were described in~\cite{Zubkov96}, \cite{Lopatin_Orel}, \cite{Lopatin_Ofree}, \cite{Lopatin_Sp}.

Given an $\NN$-graded algebra $\algA$ with the component of degree zero equal to $\FF$, denote by $\algA^{+}$ the subalgebra generated by homogeneous elements of positive degree. A set $\{a_i\} \subseteq \algA$ is a m.h.s.g.~for the unital algebra $\algA$~if and only if the $a_i$'s are $\NN$-homogeneous and $\{\ov{a_i}\}$ is a basis of $\ov{\algA}={\algA}^{+}/{(\algA^{+})^2}$. If we consider $a\in\algA$ as an element of $\ov{\algA}$, then we usually omit the bar and write $a\in\ov{\algA}$ instead of $\ov{a}$.  If for an element $a\in \algA$ we have $a=0$ in $\ov{\algA}$, then $a$ is called {\it decomposable} and we write $a\equiv0$. In other words, a decomposable element is equal to a polynomial in elements of strictly lower degree. Therefore  the highest degree $D(\algA)$ of indecomposable elements of  $\algA$  is equal to the least upper bound for the degrees of elements of a m.h.s.g.~for $\algA$.

A m.h.s.g.~for matrix $GL(2)$ and $O(2)$-invariants are known (see~\cite{Sibirskii_1968}, \cite{Procesi_1984}, \cite{DKZ_2002}). A m.h.s.g.~for matrix $GL(3)$-invariants was established in~\cite{Lopatin_Comm1} and~\cite{Lopatin_Comm2}. A m.h.s.g.~for $R^{GL(4)}$ in case $\Char{\FF}=0$ and $d=2$ was established in~\cite{Drensky_Sadikova_4x4} (see also~\cite{Teranishi_1986}). M.h.s.g.-s~for $O(n)$- and $SO(n)$-invariants of one traceless matrix over an algebraically closed field of characteristic zero were given in~\cite{Djokovic_O4} and~\cite{Djokovic_Osmall} in case $n\leq 5$.

Let us remark, that the case of characteristic zero is drastically different from the case of a ``small''{} positive characteristic. Namely,  for $0<\Char{\FF}\leq n$ we have $D(R^{GL(n)})\to\infty$ as $d\to\infty$ (see~\cite{DKZ_2002}), but in case $\Char\FF=0$ or $\Char{\FF}>n$ this statement does not hold. In the recent preprint~\cite{Lopatin_tends} it was shown that $D(R^{O(n)})\to\infty$ as $d\to\infty$ for  $0<\Char{\FF}\leq n$. In~\cite{Lopatin_O3} some estimations on $D(R^{O(3)})$ were given.

\section{Finite set of generating identities for $I(\Q,\n,\i)$}

In this section $\Q$, $\n$ and $\i$ are the same as above.
Consider some notions concerning all products of arrows of $\Q^{\D}$.

\begin{enumerate}
\item[$\bullet$] Let $\YQ$ be the semigroup (without unity) freely generated by arrows  $a\in\Q^{\D}$, where $a\in\Q_1^{\D}$, and $\YQ^{\#}=\YQ\sqcup\{1\}$.

\item[$\bullet$] Introduce some lexicographical linear order on $\YQ$ with $ab>a$ for $a,b\in\YQ$.

\item[$\bullet$] Introduce the involution ${}^T$ on $\YQ$ as follows. We set $(a)^T=a^T$, $(a^T)^T=a$ and $(a_1\cdots a_s)^T=a_s^T\cdots a_1^T\in\Y$ for all $a, a_1,\ldots,a_s\in\Q^{\D}_1$. 

\item[$\bullet$] We say that $a,b\in\Y$ are {\it cyclic equivalent} and write $a\stackrel{c}{\sim} b$ if $a=a_1a_2$ and $b=a_2a_1$ for some $a_1,a_2\in\Y^{\#}$. If $a\stackrel{c}{\sim} b$ or $a\stackrel{c}{\sim} b^T$, then we say that $a$ and $b$ are {\it equivalent} and write $a\sim b$. 


\item[$\bullet$] $\LA\Q\RA\subset \YQ$ is the set of all (non-empty) closed paths in $\Q^{\D}$ and $\LA\Q\RA^{\#}=\LA\Q\RA \sqcup\{1\}$.


\end{enumerate}%

\noindent{}We can define the so-called {\it small free}  algebra for $I(\Q,\n,\i)$:
 
\begin{enumerate}
\item[$\bullet$] Let $\FF\LA\Q\RA$ and $\FF\LA\Q\RA^{\#}$ be the free associative algebras (without and with unity, respectively) with the $\FF$-bases $\LA\Q\RA$ and $\LA\Q\RA^{\#}$, respectively. Note that elements of $\FF\LA\Q\RA$ and $\FF\LA\Q\RA^{\#}$ are {\it finite} linear combinations of monomials from $\LA\Q\RA$ and $\LA\Q\RA^{\#}$, respectively.

\item[$\bullet$] Let $u,v\in\Q_0$ be vertices. We say that $a\in\FF\LA\Q\RA$ goes from $u$ to $v$ 
if $a=\sum_i\al_i a_i$, where $\al_i\in\FF$ and $a_i$ is a path in $\Q^{\D}$ from $u$ to $v$. If $a$ goes from $u$ to $u$, then we say that $a$ is {\it
incident} to $u$. Similarly we define that a vector $\un{a}=(a_1,\ldots,a_r)$ of paths in $\Q^{\D}$  goes from $u$ to $v$ (is incident to $v$, respectively). 

\item[$\bullet$] A triple $(a,b,c)$ of paths is {\it admissible} in a vertex $v$ of $\Q$ if $a$ is incident to $v$, $b$ goes from $\i(v)$ to $v$ and $c$ goes from $v$ to $\i(v)$.  Similarly we define that a triple $(\un{a};\un{b};\un{c})=(a_1,\ldots,a_r;b_1,\ldots,b_s;c_1,\ldots,c_t)$ of vectors of paths in $\Q^{\D}$  is admissible in $v$. 

\item[$\bullet$] Define a homomorphism of algebras $\phi_{\n}:\FF\LA\Q\RA^{\#}\to \alg_{\FF}\{E,X_a\, |\,a\in\Q_1^{\D}\}$ by $1\to E$ and $a\to X_a$ for all $a\in\Q_1^{\D}$. Note that $X_{a^T}=(X_a)^T$.

\item[$\bullet$] Let $\AlgSmallQ{\n}$ be a ring with unity of commutative polynomials over $\FF$ freely generated by ``symbolic'' elements $\si_t(a)$, where $a\in\FF\LA\Q\RA$ is incident to some vertex $v$ of $\Q$,  $a$ has the coefficient $1$ in the highest term with respect to the introduced lexicographical order on $\LA\Q\RA\subset\YQ$ and $1\leq t\leq \n_v$. Define 
$$\si_t(\al a)=\al^t\si_t(a)$$
for $\al\in\FF$ and denote $\si_0(a)=1$, $\tr(a)=\si_1(a)$. Note that $\si_t(0)=0$.
\end{enumerate}

Consider the surjective homomorphism 
$$\PhiSmallQ{\n}:\AlgSmallQ{\n}\to I(\Q,\n,\i)$$ 
such that $\si_t(a) \to \si_t(\phi_n(a))$, where $a\in\FF\LA\Q\RA$ is incident to a vertex $v\in\Q$ and $1\leq t\leq \n_v$. Since
$$\si_t(\al A)=\al^t\si_t(A)$$
holds for an arbitrary $n\times n$ matrix $A$ over a commutative $\FF$-algebra and $1\leq t\leq n$, the homomorphism $\PhiSmallQ{\n}$ is well-defined. Its kernel $\KSmallQ{\n}$ is the ideal of relations for $I(\Q,\n,\i)$ in the small free algebra.


The definitions of elements $F_t(a,b)$, $P_{t,l}(a)$ and $\si_{\un{t};\un{r};\un{s}}(\un{a};\un{b};\un{c})$ from the above given Theorem~\ref{theo1} can be found in~\cite{Lopatin_JPAA2013}. In case $\un{t}=(t)$, $\un{r}=(r)$, $\un{s}=(s)$ we have $r=s$, $\un{a}=(a)$, $\un{b}=(b)$, $\un{c}=(c)$ and we write $\si_{t,r}(a,b,c)$  for $\si_{\un{t};\un{r};\un{s}}(\un{a};\un{b};\un{c})$.  Note that $\si_{\un{t};\un{r};\un{s}}(\un{a};\un{b};\un{c})$ is the {\it partial linearization} of $\si_{t,r}(a,b,c)$ in case $t=|\un{t}|$ and $r=|\un{r}|=|\un{s}|$.

\begin{theo}\label{theo1} \mylabel{theo1}
The ideal of relations $\KSmallQ{\n}$ for $I(\Q,\n,\i)\simeq \AlgSmallQ{\n} / \KSmallQ{\n}$ is generated by 
\begin{enumerate}
\item[(a)] $\si_t(a+b)=F_t(a,b)$ for $1\leq t\leq n$, where $a,b\in\FF\LA\Q\RA$ are incident to one and the same vertex;  
\item[(b)] $\si_t(a^l)=P_{t,l}(a)$ for $1\leq t\leq n$, $1<l\leq n$, where $a\in\LA\Q\RA$;

\item[(c)] $\si_t(ab)=\si_t(ba)$ for $1\leq t\leq n$, where $ab\in\LA\Q\RA$;

\item[(d)] $\si_t(a)=\si_t(a^T)$ for $1\leq t\leq n$, where $a\in\LA\Q\RA$;

\item[(e)] $\si_{\un{t};\un{r};\un{s}}(\un{a};\un{b};\un{c})=0$ for $n<|\un{t}|+2|\un{r}|\leq 2n$, where $\un{t}\in\NN_0^u$, $\un{r}\in\NN_0^v$, $\un{s}\in\NN_0^w$ ($u,v,w>0$) satisfy $|\un{r}|=|\un{s}|$, $a_i,b_j,c_k\in\LA\Q\RA$ for all $i,j,k$ and $(\un{a};\un{b};\un{c})$ is an admissible triple of vectors.
 \end{enumerate}
\end{theo}
\bigskip

Theorem~\ref{theo1} implies that the ideal $\KSmallQ{\n}$ is finitely based as T-ideal (see~\cite{Lopatin_JPAA2013} for details). 

\begin{remark}\label{rem-to-theo1} \mylabel{rem-to-theo1}
In the formulation of Theorem~\ref{theo1} we can assume that $\un{t},\un{r},\un{s}$ from relation~(e) satisfy the following conditions:
\begin{eq}\label{eq_cond1}
t_1\geq \cdots \geq t_u,\; r_1\geq \cdots \geq r_v,\; s_1\geq \cdots \geq s_w
\end{eq}
\vspace{-0.5cm}%
\begin{eq}\label{eq_cond2}
k_1,\ldots,k_q\in\{1,p,p^2,p^3,\ldots\},
\end{eq}
\vspace{-0.5cm}%
\begin{eq}\label{eq_cond3}
\text{either } |\un{k}|=n+1, \text{ or } n+1<|\un{k}|\leq 2n \text{ and } |\un{k}|-\min\{k_i\}\leq n. %
\end{eq}%
where $\un{k}=(k_1,\ldots,k_q)$ is the result of elimination from  $(t_1,\ldots,t_u,r_1,\ldots,r_v,s_1,\ldots,s_w)$ of all zero entries.  In particular, if the second case from~(\ref{eq_cond3}) holds, then $k_i\neq 1$ for all $i$.

\end{remark}
\bigskip

\begin{proof}
This statement is an analogue of Theorem 5.5 from~\cite{Lopatin_JPAA2013}, which describes the finite set of generators of the ideal of identities of the algebra $R^{O(n)}$ of matrix $O(n)$-invariants in case $p\neq2$. Earlier in this paper we have  defined $\LA\Q\RA$, $\FF\LA\Q\RA$, $\AlgSmallQ{\n}$, $\PhiSmallQ{\n}$, $\KSmallQ{\n}$ similarly to the notions $\Y$, $\FF\Y$, $\AlgSmallY{n}$, $\PhiSmallY{n}$,  $\KSmallY{n}$, respectively, from~\cite{Lopatin_JPAA2013}. Then we can define $\si\LA\Q\RA$, $\AlgLargeQ$, $\piSmallQ{\n}$, $\piLargeQ$, $\PhiAbsQ{\n}$, $\PhiLargeQ{\n}$, $\KLargeQ{\n}$ similarly to $\si\Y$, $\AlgLargeY$, $\piSmallY{n}$, $\piLargeY$, $\PhiAbsY{n}$, $\PhiLargeY{n}$, $\KLargeY{n}$, respectively. The following diagram is commutative:
$$ 
\begin{picture}(0,120)
\put(0,95){%
\put(0,-2){\vector(0,-1){58}}%
\put(15,0){\vector(3,-2){35}}%
\put(-15,0){\vector(-3,-2){35}}%
\put(-11,5){$\si\LA\Q\RA$}%
\put(-75,-33){$\AlgSmallQ{\n}$}%
\put(50,-33){$\AlgLargeQ$}%
\put(50,-40){\vector(-3,-2){35}}%
\put(-50,-40){\vector(3,-2){35}}%
\put(110,0){\vector(-3,-2){35}}%
\put(-110,0){\vector(3,-2){35}}%
\put(-24,-75){$I(\Q,\n,\i)$}%
\put(-125,5){$\KSmallQ{\n}$}%
\put(115,5){$\KLargeQ{\n}$}%
\put(3,-33){$\scriptstyle\PhiAbsQ{\n}$}%
\put(-40,-8){$\scriptstyle\piSmallQ{\n}$}%
\put(33,-8){$\scriptstyle\piLargeQ$}%
\put(-35,-48){$\scriptstyle\PhiSmallQ{\n}$}%
\put(25,-48){$\scriptstyle\PhiLargeQ{\n}$}%
\put(-20,-90){\text{Diagram 1.}}%
}%
\end{picture}
$$%
In this diagram all maps are surjections and $\KLargeQ{\n}$ is the kernel of $\PhiLarge{\n}$. The description of the ideal $\KLargeQ{\n}$ modulo free relations of $I(\Q,\n,\i)$ was obtained in~\cite{ZubkovII}. In Lemma~6.1 of~\cite{Lopatin_Ofree} it was shown that do not exist non-trivial free relations for $I(\Q,\n,\i)$ and therefore the ideal $\KLargeQ{\n}$ is generated by
\begin{eq}\label{eq_KLarge}\mylabel{eq-KLarge} 
\si_{t,r}(a,b,c),
\end{eq}%
where $(a,b,c)$ is an admissible triple of elements of $\FF\LA\Q\RA$ in some vertex $v\in\Q_0$ and $t+2r>\n_v$, $t,r\geq0$ (see Theorem~6.2 of~\cite{Lopatin_Ofree}). Note that the given set of identities is not finite. We start the proof with this description of the ideal $\KLargeQ{\n}$ in the same way as the proof in~\cite{Lopatin_JPAA2013} started with Theorem~6.1.


Given $l\geq0$, we denote by $J_{l}$ the ideal of $\AlgLargeQ$ generated by $\si_{t,r}(a,b,c)$, where $l=t+2r$ and  $(a,b,c)$ is an admissible triple of some elements of $\FF\LA\Q\RA$.

\medskip
\noindent{}{\bf Fact 1.} {\it Given an admissible triple $(x_0,x;y;z)$ of elements of $\LA\Q\RA$ and $k,t,r\geq0$, we have $\si_{k,t;r;r}(x_0,x;y;z)\in J_{t+2r}$.
}
\medskip

\noindent{}Fact 1 is an analogue of Lemma~6.9 of~\cite{Lopatin_Ofree}. Its proof is based on the following formula
\begin{eq}\label{eq-fact1} 
\begin{array}{c}
\si_{k,t;r;r}(x_0,x;y;z) \\
=\sum\limits(-1)^{\al_0+k}\si_{\al_0}(x_0)\, \si_{\Delta'}(x,x_0^i x; y, x_0^i y, y (x_0^T)^i, x_0^i y (x_0^T)^j; z)_{i,j>0},\\
\end{array}
\end{eq}%
where the sum ranges over the certain list of multidegrees $\Delta$. Note that $(x,x_0^i x; y, x_0^i y, y (x_0^T)^i, x_0^i y (x_0^T)^j; z)$ is an admissible triple. Then we can just repeat the proof.

\medskip
\noindent{}{\bf Fact 2.} {\it Given $x,y_0,y,z\in\Y$ and $t,r,s\geq0$, we have $\si_{t;r,s;r+s}(x;y_0,y;z)\in J_{t+r+2s}$.}
\medskip

\noindent{}Fact 2 is an analogue of Lemma~6.12 of~\cite{Lopatin_Ofree}. Its proof is based on the following formula
\begin{eq}\label{eq-fact2}
\si_{t;r,s;r+s}(x;y_0,y;z) = \sum\limits (-1)^{\al_2+r} \si_{\Delta}(x, y_0 z, y_0 z^T; y, y_0x^T;  z),
\end{eq}%
where the sum ranges over the certain list of multidegrees $\Delta$ with $|\Delta|=t+r+2s$. Note that $(x, y_0 z, y_0 z^T; y, y_0x^T;  z)$ is an admissible triple. Then we can just repeat the proof.

\medskip
\noindent{}{\bf Fact 3.} {\it Given $\un{t}\in\NN^u_0$, $\un{r}\in\NN^v_0$, $\un{s}\in\NN^w_0$ with $t=|\un{t}|$ and $r=|\un{r}|=|\un{s}|$ ($u,v,w>0$). Then for $f=\si_{\un{t};\un{r};\un{s}}(\un{a};\un{b};\un{c})$ from $\AlgLargeQ$, where $(\un{a};\un{b};\un{c})$ is an admissible triple of vectors from $\LA\Q\RA$, we have
 \begin{enumerate}
  \item[(a)] $f\in J_{t+2r-t_1}$;
  
  \item[(b)] $f\in J_{t+2r-r_1}$ and $f\in J_{t+2r-s_1}$.
 \end{enumerate}
}
\medskip

\noindent{}Fact 3 is an analogue of Lemma~6.14 of~\cite{Lopatin_Ofree}. Its proof follows from Fact 1 and 2 in the same way as in~\cite{Lopatin_Ofree}.

Fact~3 implies that the ideal of relations $\KLargeQ{\n}$ for $I(\Q,\n,\i)\simeq \AlgLargeQ / \KLargeQ{\n}$ is generated by 
\begin{enumerate}
\item[$\bullet$] $\si_{t,r}(a,b,c)=0$ for $n<t+2r\leq 2n$, $t,r\geq0$, and $(a,b,c)$ is an admissible triple of elements of $\FF\LA\Q\RA$; 

\item[$\bullet$] $\si_t(b)=0$, where $t>2n$ and $b\in\EQ$ is a closed path in $\Q^{\D}$. 
\end{enumerate}

The rest of the proof the same as in~\cite{Lopatin_Ofree}. Note that the restriction that $p\neq2$ is only required in the formulation of Theorem~6.1 of~\cite{Lopatin_Ofree} and does not play any role in the rest of the proof.
\end{proof}

\section{The case of dimension vector $(2,\ldots,2)$}

In this section $\Q$, $\n$ and $\i$ are the same as in Section~1.2. Moreover we assume that $\n=(2,\ldots,2)$. For short we write $\n=2^l$, where $l$ is the number of vertices in $\Q$. Denote by $\tr\LA Q\RA$ the free associative commutative algebra generated by free variables $\tr(a)$, where $a\in\LA\Q\RA$ is a closed path in $\Q^\D$. As we know from Section 1.2, $I(\Q,2^l,\i)$ is generated by $\tr(X_a)$, $\det(X_b)$, where $a,b\in\LA\Q\RA$ and $b$ is without self-intersections. In case $\sum_{i=1}^r\al_i\tr(X_{a_i})\equiv 0$ for some incident elements $a_1,\ldots,a_r$ of $\LA\Q\RA$ we write $\sum_i\tr(a_i)\equiv 0$.

For $a\in\FF\LA\Q\RA$ we denote 
$$\ov{a}=a-a^T.$$
For short we also write $\tr(a\ov{b}c)$ for  $\tr(abc)-\tr(ab^Tc)\in\tr\LA\Q\RA$.

\begin{theo}\label{theo-222} \mylabel{theo-222} 
Assume that $h=\sum_i\al_i a_i$ is an incident element of $\FF\LA\Q\RA$, where $a\in\LA\Q\RA$ and $\al_i\in\FF$. Then $\tr(X_a)\equiv 0$ in $I(\Q,2^l,\i)$ if and only if $a$ is a linear combination of the following elements of $\FF\LA\Q\RA$:
\begin{enumerate}
\item[$(1)$] $\tr(uv)\equiv \tr(vu)$, where $uv\in\LA\Q\RA$;

\item[$(2)$] $\tr(a_{\si(1)}\ldots a_{\si(t)}) \equiv \tr(a_1\cdots a_t)$, where $t\geq2$ and $\si\in S_t$;

\item[$(3)$] $\tr(a_1^2 a_2)\equiv0$;

\item[$(4)$] if $p=2$, then $\tr(a^2)\equiv0$,\\
if $p\neq2$, then $\tr(a_1\cdots a_4)\equiv0$.

\item[$(5)$] $\tr(a\ov{b}\ov{c})\equiv0$;

\item[$(6)$] $\tr(a)\equiv \tr(a^T)$. 
\end{enumerate}
Here $a,a_1,\ldots,a_4\in\LA\Q\RA$ are incident and $(a,b,c)$ is an admissible triple of elements of $\LA\Q\RA$.
\end{theo}

\begin{proof} We consider $\tr\LA\Q\RA$ as a subalgebra of $\AlgSmallQ{2^l}$ and we work in $\AlgSmallQ{2^l}$. 
We have 
$$\si_{1,1}(a,b,c) = \tr(a\ov{b}\ov{c}) - \tr(a)\tr(b\ov{c}) \text{ and }$$
$$\si_{(2,1)}(a_1,a_2)=\tr(a_1^2 a_2) - \tr(a_1a_2)\tr(a_1) + \det(a_1)\tr(a_2).$$
Then we obtain that (3) and (5) holds. It is well known (and it is easy to verify), that (3) implies (4) in case $p\neq2$. Equivalence (4) in case $p=2$ follows from $$\tr(a_1^2)=P_{1,2}(a_1)=\tr(a)^2 - 2 \si_2(a).$$
The case by case consideration of the finite set of identities from Theorem~\ref{theo1} completes the proof.
\end{proof}

\end{document}